\documentclass{amsart}[12pt]
\usepackage{amsmath,amssymb,amscd,latexsym}
\usepackage[all]{xy}

\newcommand{\C}{{\mathbb{C}}}

\newcommand{\N}{{\mathbb{N}}}

\newcommand{\R}{{\mathbb{R}}}
\newcommand{\T}{{\mathbb{T}}}
\newcommand{\Z}{{\mathbb{Z}}}

\newcommand{\Ch}{{\mathcal C}}

\newcommand{\Oh}{{\mathcal O}}

\newcommand{\Zh}{{\mathcal Z}}

\newcommand{\Aff}{\mathrm{Aff}}

\newcommand{\ev}{\mathrm{ev }}

\newcounter{number}[section]

\newenvironment{nummer}{\refstepcounter{number}{\noindent\arabic{section}.\arabic{number}}}{}

\newcommand{\bn}{\noindent \begin{nummer} \rm}
\newcommand{\en}{\end{nummer}}

\newenvironment{thms}{\noindent {\sc Theorem:} \it}{}

\newenvironment{props}{\noindent {\sc Proposition:} \it}{}
\newenvironment{dfs}{\noindent {\sc Definition:} \it}{}
\newenvironment{cors}{\noindent {\sc Corollary:} \it}{}

\newenvironment{qus}{\noindent {\sc Question:} }{}
\newenvironment{nproof}{\noindent {\sc Proof:}}{\mbox{}\hfill 
\rule[-.2ex]{.25em}{1.8ex}}

\begin{document}

\title[$\Zh$-stable ASH Algebras]{{\sc $\Zh$-stable ASH Algebras}}

\author{Andrew S. Toms}
\address{Department of Mathematics and Statistics, University of New Brunswick\\
Frederic- 
\indent ton, New Brunswick\\
E3B 5A3}

\email{atoms@math.unb.ca}

\author{Wilhelm Winter}
\address{Mathematisches Institut der Universit\"at M\"unster\\
Einsteinstr. 62\\ D-48149 M\"unster}

\email{wwinter@math.uni-muenster.de}

\date{August 12, 2005}
\subjclass{46L85, 46L35}
\keywords{nuclear $C^*$-algebras, K-theory,  
classification}
\thanks{{\it Supported by:} DFG (through the SFB 478), EU-Network  Quantum Spaces - Noncommutative 
\indent Geometry (Contract No. HPRN-CT-2002-00280) and  NSERC Postdoctoral Fellowship}

\setcounter{section}{-1}

\begin{abstract}
The Jiang--Su algebra $\mathcal{Z}$ has come to prominence in the classification
program for nuclear $C^*$-algebras of late, due primarily to the fact that Elliott's 
classification conjecture predicts that all simple, separable, and nuclear $C^*$-algebras
with unperforated $\mathrm{K}$-theory will absorb $\mathcal{Z}$ tensorially (i.e., will be $\mathcal{Z}$-stable).
There exist counterexamples which suggest that the conjecture
will only hold for simple, nuclear, separable and $\mathcal{Z}$-stable $C^*$-algebras.
We prove that virtually all classes of nuclear $C^*$-algebras for which the Elliott
conjecture has been confirmed so far, consist of $\mathcal{Z}$-stable $C^*$-algebras.
This result follows in large part from the following theorem, also proved herein:
separable and approximately divisible $C^*$-algebras are $\mathcal{Z}$-stable.
\end{abstract}

\maketitle

\section{Introduction}

The Jiang--Su algebra $\mathcal{Z}$ is a simple, separable, unital and nuclear $C^*$-algebra
$\mathrm{KK}$-equivalent to $\mathbb{C}$ (\cite{JS1}).  Since its discovery in 1995 there has been a steady
accumulation of evidence linking $\mathcal{Z}$ to Elliott's program to classify separable, nuclear
$C^*$-algebras via $\mathrm{K}$-theoretic invariants:  in \cite{JS1}, Jiang and Su prove 
that simple, infinite-dimensional AF algebras and Kirchberg algebras (simple, nuclear, 
purely infinite and satisfying the Universal Coefficients Theorem) are $\mathcal{Z}$-stable, 
i.e., for any such algebra $A$ one has an isomorphism $\alpha:A \to A \otimes \mathcal{Z}$. 
This in particular implies that the purely infinite algebras covered by Kirchberg's classification 
of $\Oh_{2}$-stable $C^{*}$-algebras and by Kirchberg--Phillips classification are $\Zh$-stable, cf.\ \cite{K}.   
In \cite{GJS}, Gong, Jiang, and Su showed that if a simple, unital and nuclear $C^*$-algebra $A$ has
a weakly unperforated ordered $\mathrm{K}_0$-group, then the ordered $\mathrm{K}_0$-groups
of $A$ and $A \otimes \mathcal{Z}$ are isomorphic.  All known counterexamples to Elliott's 
classification conjecture fail to be $\mathcal{Z}$-stable.

In the present paper --- a natural sequel to \cite{TW} --- we prove that as of this moment,
all classes of non-type-I $C^*$-algebras for which the Elliott conjecture is confirmed consist entirely of $\mathcal{Z}$-stable
algebras.  In the approximately homogeneous (AH) case, this result follows for the most part from \cite{EGL}, \cite{EGL2} and our Theorem 3.1, which states that separable
and approximately divisible $C^*$-algebras are $\mathcal{Z}$-stable.  One must then expend considerable effort
in proving that certain approximately subhomogeneous (ASH) $C^*$-algebras which are not approximately 
homogeneous are nevertheless $\mathcal{Z}$-stable --- such algebras need not be approximately divisible.  
To the best of the authors' knowledge, the classification results of \cite{JS1}, \cite{JS2}, \cite{EJS},  \cite{T}, \cite{M}, \cite{R}, \cite{I}, 
\cite{Ts} exhaust the known classes of properly ASH (i.e., \emph{not} AH) algebras 
for which the Elliott conjecture is confirmed.  We prove that some of these
classification results in fact cover approximately divisible $C^*$-algebras, and so are $\mathcal{Z}$-stable
by Theorem 3.1.  For the remaining classes of ASH algebras, we develop an approach to proving $\mathcal{Z}$-stability
which should remain applicable as more general classification results for ASH algebras arise.   

We feel that our results, combined with the counterexamples of \cite{R3}, \cite{To1} and 
\cite{To2}, demonstrate the necessity for stabilisation by $\mathcal{Z}$ in the classification
program for separable, nuclear $C^*$-algebras.

\emph{Acknowledgements:} The authors would like to thank George Elliott and Mikael R{\o}rdam for
many inspiring conversations, and, in particular, for the suggestion that separable and approximately
divisible $C^*$-algebras might be $\mathcal{Z}$-stable, regardless of whether or not they are nuclear.

\section{$\mathcal{Z}$-stability and an augmented invariant}

In this section we review the effect of tensoring with $\mathcal{Z}$
on the Elliott invariant of a simple, unital and nuclear $C^*$-algebra $A$, and
examine the effect of this operation on an augmented version of the invariant.
Our conclusion, predictably, is that tensoring with $\mathcal{Z}$ has no
effect whatever provided that (the invariant of) $A$ is sufficiently well-behaved.

\bn
Let $A$ be a simple and nuclear $C^*$-algebra.  Define an invariant
\[
\mathrm{I}(A) := \left( (\mathrm{K}_*A,\mathrm{K}_*A^+),\mathrm{T}^+A,r_A,\mu_A \right),
\]
where $(\mathrm{K}_*A,\mathrm{K}_*A^+)$ is the (pre-)ordered topological $\mathrm{K}$-theory
of $A$ (this includes the order on the direct sum $\mathrm{K}_*A := \mathrm{K}_0A \oplus \mathrm{K}_1A$
coming from partial unitaries, which we review below),  
$\mathrm{T}^+A$ is the (possibly empty) space of tracial functionals,
$r_A$ is the pairing between $\mathrm{T}^+A$ and $\mathrm{K}_0A$ given by evaluation, and
$\mu_A:\mathrm{T}^+A \to \mathbb{R}^+$ is the trace-norm map defined by 
\[
\mu_A(\tau)= \mathrm{sup}_{a \in \mathcal{B}_1(A)} \tau(a) \, ,
\]
where $\mathcal{B}_1(A)$ is the unit ball of $A$.  If $A$ is unital, then we 
replace $\mathrm{T}^+A$ with the space of tracial states $\mathrm{T}A$, and 
ignore $\mu_A$ (it is identically equal to one).  We include instead the 
$\mathrm{K}_0$-class of the unit, $[\mathbf{1}_A]$.  The space $\mathrm{T}^{+}A$ 
and the pairing $r_A$ are only relevant in the setting of stably finite $C^*$-algebras.

The invariant $\mathrm{I}(A)$ is a (slightly) augmented version of the usual Elliott invariant.
The latter, say $\mathrm{Ell}(A)$, is obtained from $\mathrm{I}(A)$ by considering only the
order structure on $\mathrm{K}_0A$, rather than $\mathrm{K}_*A$. 

To prepare the next proposition, we review the definition of 
$\mathrm{K}_*A^+$.  An element $u$ of a unital $C^*$-algebra $A$ is called a
\emph{partial unitary} if it is both normal and a partial isometry.  For such 
a $u$, the element $u + (\mathbf{1}_A - u u^*)$ is a unitary in $A$.
Say that an element 
\[
x \oplus y \in \mathrm{K}_0A \oplus \mathrm{K}_1A = \mathrm{K}_*A
\]
is
positive if
\begin{enumerate}
\item $x = [p]$ for some projection $p \in \mathrm{M}_n(A)$, and
\item there is a partial unitary $v \in \mathrm{M}_k(A)$ such that
\[
[vv^*] \leq [p]  \ \ \mathrm{and} \ \ [v + (\mathbf{1}_{\mathrm{M}_k(A)} - vv^{*})]_1 =y.
\]
\end{enumerate}
\en

\bn
\begin{props}\label{Z-stable-invariant}
Let $A$ be a simple, separable, unital, and nuclear $C^*$-algebra.  Then $\mathrm{I}(A) \cong
\mathrm{I}(A \otimes \mathcal{Z})$ iff the (pre-)ordered group $(\mathrm{K}_0A, \mathrm{K}_0A^+)$
is weakly unperforated and the order structure on $\mathrm{K}_*A$ is the strict order
coming from the direct summand $\mathrm{K}_0A$ of $\mathrm{K}_*A$.
\end{props}

\begin{nproof}
In \cite{GJS} it is shown that $\mathrm{Ell}(A) \cong \mathrm{Ell}(A \otimes \mathcal{Z})$ whenever
$(\mathrm{K}_0A,\mathrm{K}_0A^+)$ is weakly unperforated, so
that any difference between $\mathrm{I}(A)$ and $\mathrm{I}(A \otimes \mathcal{Z})$ must occur 
in the cone $\mathrm{K}_*A^+$.

Theorem 2 of \cite{J} states that the natural map 
\[
\eta: \mathcal{U}(B \otimes \mathcal{Z})/\mathcal{U}(B \otimes \mathcal{Z})_0 \to \mathrm{K}_1(B \otimes \mathcal{Z})
\]
is an isomorphism
whenever $B$ is a unital $C^*$-algebra.  For a projection $p \in B$, the map
\[
\iota:\mathrm{K}_1 (pBp) \to \mathrm{K}_1 B
\]
given by sending the class of a unitary $v \in \mathrm{M}_k(pBp)$ to the class of the unitary
\[
v + (\mathbf{1}_{\mathrm{M}_k(B)} - \mathbf{1}_{\mathrm{M}_k(pBp)}) \in \mathrm{M}_k(B)
\]
is surjective whenever $p$ is full.  Let $q \in A \otimes \mathcal{Z}$ be a projection, and 
let $y \in \mathrm{K}_1(A \otimes \mathcal{Z})$.  Since $A$ is simple, $q$ is full and the map 
\[
\iota:\mathrm{K}_1 (q(A \otimes \mathcal{Z})q) \to \mathrm{K}_1 (A \otimes \mathcal{Z})
\]
is surjective, hence there is a unitary $v \in \mathrm{M}_k(q(A \otimes \mathcal{Z})q)$, some $k \in
\mathbb{N}$, such that $\iota([v]_{1}) = y$.  By Lemma 3.2 of \cite{J}, $\mathcal{Z}$-stability
passes to corners (in fact, to hereditary subalgebras --- see section 3 
of \cite{TW}), whence the map
\[
\eta: \mathcal{U}(q(A \otimes \mathcal{Z})q)/\mathcal{U}(q(A \otimes \mathcal{Z})q)_0
\to \mathrm{K}_1(q(A \otimes \mathcal{Z})q)
\]
is an isomorphism.  It follows that $v$ may be chosen to lie in $q(A \otimes \mathcal{Z})q$.
Thus, $([q],y)$ is positive for every $y$, and $\mathrm{K}_* (A \otimes \mathcal{Z})$ 
has the strict order coming from $\mathrm{K}_0 (A\otimes \mathcal{Z})$.
\end{nproof}
\en

Propsition \ref{Z-stable-invariant} does not address non-unital algebras, and cannot
be adapted immediately to a non-unital algebra with unit adjoined.  Nevertheless,
we can say something about $\mathrm{I}(A \otimes \mathcal{Z})$ in this case:  
if $\gamma$ is the unique normalised trace on $\mathcal{Z}$, then 
the isomorphism $\zeta: \mathrm{T}^+A \to \mathrm{T}^+ (A \otimes \mathcal{Z})$ given by
$\zeta(\tau) = \tau \otimes \gamma$ preserves the trace-norm map.  

\section{Approximate divisibility}

In this section we show that separable, approximately divisible $C^*$-algebras are $\mathcal{Z}$-stable, 
answering a question posed by Jiang in \cite{J}. We first recall the definition of approximate divisibility, 
then give a new version of Theorem 2.2 of \cite{TW} to make  the special 
inductive limit decomposition of $\mathcal{Z}$ available for our purposes. Following the terminology of 
\cite{TW}, we denote by $\mathcal{Q}(A)$ the quotient $\prod_{\N}A/\bigoplus_{\N}A$ for any $C^{*}$-algebra $A$; 
$\mathcal{M}(A)$ will be the multiplier algebra of $A$.

\bn
\begin{dfs}
A $C^{*}$-algebra $A$ is said to be approximately divisible, if, for any $N \in \mathbb{N}$, there 
is a sequence of unital $*$-homomorphisms $\mu_{n}: \mathrm{M}_{N} \oplus \mathrm{M}_{N+1} \to \mathcal{M}(A)$ which 
is approximately central for $A$, i.e., $\|[\mu_{n}(x),a]\| \stackrel{n \to \infty}{\longrightarrow} 
0$ for all $a \in A$ and $x \in \mathrm{M}_{N} \oplus \mathrm{M}_{N+1}$.
\end{dfs}
\en

\bn
\begin{props}\label{inductive-limit-intertwining}
Let $A$ and $\mathcal{D}$ be  separable  $C^{*}$-algebras, $\mathcal{D}$ unital, strongly self-absorbing 
and $\mathrm{K}_{1}$-injective. Suppose that $\mathcal{D}$ can be written as the closure of an increasing union of 
nuclear $C^{*}$-algebras $\mathcal{D}_{i}, \, i \in \mathbb{N}$.  If, for each $i \in \mathbb{N}$, there exists a 
unital $*$-homomorphism $\gamma_{i}:\mathcal{D}_{i} \to \mathcal{Q}(\mathcal{M}(A)) \cap A'$, then there is an 
isomorphism $\varphi: A \to A \otimes \mathcal{D}$ and $\varphi \approx_{\mathrm{a.u.}} \mathrm{id}_{A} \otimes \mathbf{1}_{\mathcal{D}}$.
\end{props}

\begin{nproof}
By the Choi--Effros lifting theorem we can lift the $\gamma_{i}$ to u.c.p.\ maps 
\[
\textstyle
\bar{\gamma}_{i}: \mathcal{D}_{i} \to \prod_{\mathbb{N}} \mathcal{M}(A) \, .
\]
Denote the components of $\bar{\gamma}_{i}$ by $\bar{\gamma}_{i,n}, \, n \in \mathbb{N}$. Each 
$\bar{\gamma}_{i,n}$ is a nuclear u.c.p.\ map, so it can be approximated pointwise by 
finite rank u.c.p.\ maps. These in turn may be extended to u.c.p.\ maps $\bar{\gamma}_{i,n}^{(k)}: 
\mathcal{D} \to \mathcal{M}(A)$ by Arveson's extension theorem. Fix $i$ for a moment; using separability of 
$\mathcal{D}_{i}$ we can  choose a suitable subsequence $(k_{n})_{n \in \mathbb{N}}$ of $(k)_{k \in \mathbb{N}}$, 
such that the u.c.p.\ map 
\[
\textstyle
\tilde{\gamma}_{i}:= (\bar{\gamma}_{i,n}^{(k_{n})})_{n \in \mathbb{N}}: \mathcal{D} \to \prod_{\mathbb{N}} \mathcal{M}(A) \, ,
\]
when restricted to $\mathcal{D}_{i}$, also lifts $\gamma_{i}$. Again, let $\tilde{\gamma}_{i,n}$ denote the components of $\tilde{\gamma}_{i}$.\\
But now, using the properties of the $\gamma_{i}$ and the fact that the $\mathcal{D}_{i}$ exhaust all 
of $\mathcal{D}$, it is straightforward to construct subsequences $(i_{k})_{k \in \mathbb{N}}$ of $(i)_{i \in \mathbb{N}}$ 
and  $(n_{k})_{k \in \mathbb{N}}$ of $(n)_{n \in \mathbb{N}}$, such that the u.c.p.\ map 
\[
\textstyle
\tilde{\gamma}:= (\tilde{\gamma}_{i_{k},n_{k}})_{k \in \mathbb{N}}: \mathcal{D} \to  \prod_{\mathbb{N}} \mathcal{M}(A) 
\]
induces a $*$-homomorphism 
\[
\textstyle
\gamma':\mathcal{D} \to \mathcal{Q}(\mathcal{M}(A)) \cap A' \, .
\]
The assertion now holds by Theorem 7.2.2 of \cite{R1}.
\end{nproof}
\en

\bn
\label{divisible-Z-stable}
\begin{thms}
Let $A$ be a separable and approximately divisible $C^{*}$-algebra. Then there is an isomorphism 
$\varphi: A \to A \otimes \mathcal{Z}$ and $\varphi \approx_{\mathrm{a.u.}} \mathrm{id}_{A} \otimes \mathbf{1}_{\mathcal{Z}}$.
\end{thms}

\begin{nproof}
By \cite{JS1}, Proposition 2.5, $\mathcal{Z}$ may be written as an increasing union of dimension drop 
algebras $I[p_{j},p_{j}q_{j},q_{j}]$ with relatively prime integers $p_{j},q_{j}$.
For each $j \in \mathbb{N}$ there are integers $m_{j},n_{j}$ such that $1= m_{j}\cdot p_{j}+n_{j} \cdot 
q_{j}$, since  $p_{j},q_{j}$ are relatively prime. Define $N_{j}:= |m_{j}| \cdot p_{j} + |n_{j}|\cdot q_{j}$ 
and unital $*$-homomorphisms
\[
\gamma'_{j}:I[p_{j},p_{j}q_{j},q_{j}] \to \mathrm{M}_{N_{j}} \oplus \mathrm{M}_{N_{j}+1} 
\]
by 
\[
\gamma'_{j} := \textstyle(\bigoplus_{1}^{|m_{j}|} \ev_{0} \, \oplus \, \bigoplus_{1}^{|n_{j}|} 
\ev_{1}) \oplus (\bigoplus_{1}^{|m_{j}|+m_{j}} \ev_{0} \, \oplus \, \bigoplus_{1}^{|n_{j}|+n_{j}} \ev_{1}) \, .
\]
By definition of approximate divisibility, for each $j \in \mathbb{N}$ there exists a unital $*$-homomorphism 
$$\iota_{j}: \mathrm{M}_{N_{j}} \oplus \mathrm{M}_{N_{j}+1} \to \mathcal{Q}(\mathcal{M}(A)) \cap A' \, . $$ 
Define unital $*$-homomorphisms 
$$\gamma_{j}: I[p_{j},p_{j}q_{j},q_{j}] \to \mathcal{Q}(\mathcal{M}(A)) \cap A'$$
by $\gamma_{j}:= \iota_{j} \circ \gamma_{j}'$ 
and the result will follow from Proposition \ref{inductive-limit-intertwining}. 
\end{nproof}
\en

\section{AH algebras}

By \cite{EGL2}, simple AH algebras of bounded topological dimension are approximately divisible, 
so Theorem \ref{divisible-Z-stable} yields

\bn
\begin{cors}
A separable, unital and simple AH algebra of finite topological dimension is $\mathcal{Z}$-stable.
\end{cors}
 
Note that the corollary covers precisely the algebras classified in \cite{EGL}.  
There is also a number of classification results for non-unital and simple AH algebras, and for
non-simple AH algebras.  In recent work Ivanescu (\cite{I}) has proved that simple $C^*$-algebras
which are stably isomorphic to approximately interval (AI) algebras
are classified up to isomorphism by the invariant $\mathrm{I}(\bullet)$ 
described in Section 1 of the present paper, generalising work of I.\ Stevens \cite{S}.  AI algebras are approximately
divisible (cf. \cite{El4}), as are their tensor products with the compact operators
(Corollary 3.2, \cite{TW}).  It follows that the stabilization $A \otimes \mathcal{K}$ 
of an algebra $A$ as treated in \cite{I} is $\mathcal{Z}$-stable.  Another application
of Corollary 3.2 of \cite{TW} shows that $A$ itself must be $\mathcal{Z}$-stable.  \\
In the non-simple case, K.\ Stevens classifies certain non-simple
and approximately divisible AI algebras in \cite{S2}, which, by Theorem \ref{divisible-Z-stable},
are $\mathcal{Z}$-stable.  There is also the impressive result of Dadarlat and Gong (\cite{DG}), which
classifies certain ASH algebras of real rank zero.  These are frequently approximately divisible,
subject to $\mathrm{K}$-theoretic conditions which we will not describe here.
\en

\bn
The converse of Theorem \ref{divisible-Z-stable} is not true; in fact, $\mathcal{Z}$ itself is a counterexample 
($\mathcal{Z}$ is not approximately divisible, since it is unital and projectionless). However, it seems natural 
to ask for a converse at least in the case of an abundance of projections:
 
\begin{qus}
Are approximate divisibility and $\mathcal{Z}$-stability equivalent 
for simple, unital, nuclear, non-type-I $C^*$-algebras of real rank zero?
\end{qus}
\en

\bn
The idea that $\mathcal{Z}$-stability should entail classifiability via the Elliott invariant 
also suggests the following

\begin{qus}\label{AH-tensor-Z-sdg}
Let $A$ be a separable, unital and simple AH algebra.  Is $A \otimes \mathcal{Z}$
an AH algebra of finite topological dimension?
\end{qus}

As it turns out, we can make some immediate progress on Question \ref{AH-tensor-Z-sdg}.
\en

\bn
\begin{props}\label{homKgroup}
Let 
\[
A := p(\Ch(X) \otimes \mathcal{K})p
\]
be a $\mathrm{rank}(p)$-homogeneous $C^*$-algebra over
a connected compact Hausdorff space $X$.  Then,
\[
\left((\mathrm{K}_0(A \otimes \mathcal{Z}), {\mathrm{K}_0(A \otimes \mathcal{Z})}^+, 
[p \otimes \mathbf{1}_{\mathcal{Z}}]), \mathrm{K}_1(A \otimes Z) \right)
\]
is a weakly unperforated graded ordered group with the strict order
coming from $\mathrm{K}_0$.  Furthermore, the strictly positive elements
of $\mathrm{K}_0(A \otimes \mathcal{Z})$ are precisely the images under
the map $\mathrm{K}_0(\mathrm{id}_A \otimes \mathbf{1}_{\mathcal{Z}})$ of
those elements in $\mathrm{K}_0(A)$ having strictly positive virtual dimension.
\end{props}

\begin{nproof}  Since $p$ is a full projection, one may repeat the proof of Proposition
\ref{Z-stable-invariant} to conclude that $\mathrm{K}_*(A \otimes \mathcal{Z})$ has
the strict order coming from $\mathrm{K}_0$.  

The elements of $\mathrm{K}_0(A)$ can be thought of
as differences of stable isomorphism classes of complex vector bundles
over $X$.  Let $x=[\gamma]-[\omega]$ be such a difference.  The quantity
$\mathrm{rank}(\gamma)-\mathrm{rank}(\omega)$ is known as the virtual dimension 
of $x$.  If $\mathrm{rank}(\gamma)-\mathrm{rank}(\omega) > 0$, then there exists
$N \in \mathbb{N}$ such that $nx \in \mathrm{K}^0(X)^+$ for all $n \geq N$ (cf. 
Theorem 8.1.2 of \cite {H}).  If
$\mathrm{rank}(\gamma)-\mathrm{rank}(\omega) = 0$, then $mx \in \mathrm{K}^0(X)^+$
for some $m \in \mathbb{N}$ if and only if $mx = 0$.  If 
$\mathrm{rank}(\gamma)-\mathrm{rank}(\omega) < 0$, then no positive 
integral multiple of $x$ is positive in $\mathrm{K}_0(A)$.
Note that $(\mathrm{K}_0(A),\mathrm{K}_0(A)^+,[p])$ is a simple
ordered group --- every non-zero positive element is represented
by a difference $x = [\gamma]-[\omega]$ where 
$\mathrm{rank}(\gamma)-\mathrm{rank}(\omega) > 0$, and the stability 
properties of vector bundles imply that every sufficiently large
natural number multiple of $x$ will dominate a second fixed element
$y \in \mathrm{K}^0(X)$.

By Corollary 1.3 of \cite{GJS}, the inclusion $\iota: A \to A \otimes \mathcal{Z}$
given by $a \mapsto a \otimes \mathbf{1}_{\mathcal{Z}}$ induces a group isomorphism
$\iota_*:\mathrm{K}_0(A) \to \mathrm{K}_0(A \otimes \mathcal{Z})$.   
Theorem 1.4 in \cite{GJS} states that, with $\iota$ as above and $x \in
\mathrm{K}_0(A)$, one has $\iota_*(x) >0$ if and only if $nx >0$ for some
$n \in \mathbb{N}$.  The hypotheses of this theorem are that $A$ is
simple and unital, but an examination of the proof shows that these can
be weakened to the assumptions that $A$ is stably finite and that
$\mathrm{K}_0(A)$ is a simple ordered group.  Applying this new version
of Theorem 1.4 of \cite{GJS} to our situation yields the following:  
\[
(\mathrm{K}_0(A \otimes \mathcal{Z}), {\mathrm{K}_0(A \otimes \mathcal{Z})}^+, 
[p \otimes \mathbf{1}_{\mathcal{Z}}])
\]
is a simple ordered group such that $\iota_*([\gamma]-[\omega]) > 0$ 
iff $\mathrm{rank}(\gamma)-\mathrm{rank}(\omega) >0$.
\end{nproof}
\en

\bn
\label{torsion-free interpolation}
With $A$ as in the proposition above, we have that 
\[
\left((\mathrm{K}_0(A \otimes \mathcal{Z}), {\mathrm{K}_0(A \otimes \mathcal{Z})}^+, 
[p \otimes \mathbf{1}_{\mathcal{Z}}]), \mathrm{K}_1(A \otimes Z) \right)
\]
has the strict order coming from $\mathrm{K}_0(A \otimes Z)$, and 
that 
\[
(\mathrm{K}_0(A \otimes \mathcal{Z}), {\mathrm{K}_0(A \otimes \mathcal{Z})}^+, 
[p \otimes \mathbf{1}_{\mathcal{Z}}])
\]
is a finitely generated ordered abelian group order isomorphic to
\[
G := \mathbb{Z}^{n} \oplus (\mathrm{K}_{0}(A\otimes \Zh))_{\mathrm{tor}} = \mathbb{Z}^n \oplus G_{\mathrm{tor}} \, ,
\]
where an element $(x_1,\ldots,x_n) \oplus g \in G$ is strictly positive if and only if
$x_1 >0$; the $x_1$ co-ordinate is the virtual dimension of a $\mathrm{K}_0$-class.  
Let $(G,G^+)$ denote this ordered group, and note that in general it does not have the Riesz interpolation
property.  It is true, however, that if $a_1,a_2 \leq b_1,b_2$ in $(G,G^+)$ and 
if the first free co-ordinates of $b_1$ and $b_2$ exceed those
of both $a_1$ and $a_2$ by at least two, then there is an interpolating element
$g$ (i.e., $a_1,a_2 \leq g \leq b_1,b_2$) --- any $g$ with first free co-ordinate
strictly greater than the first free co-ordinates of the $a_i$'s and strictly less
than the first free co-ordinates of the $b_i$'s will serve.  In fact, the same is 
true when $G \oplus \mathrm{K}_1(A \otimes \mathcal{Z})$ is equipped with the 
strict order from $(G,G^+)$.  We use these observations to prove the main results
of this section.  We recall first some standard notions and notation.  Let $\Aff(\mathrm{T}(A))$ 
denote the space of continuous real-valued affine functions on $\mathrm{T}(A)$,
and let $W(A)$ denote the Cuntz monoid of generalised Murray--von Neumann equivalence
classes of positive elements in the algebraic direct limit  
\[
\mathrm{M}_{\infty}(A) = \lim_{i \to \infty}(\mathrm{M}_i(A),\psi_i) \, ,
\]
where $\psi_i:\mathrm{M}_i(A) \to \mathrm{M}_{i+1}(A)$ denotes inclusion as the
upper left corner.
One says that $W(A)$ is almost unperforated if $x \leq y$ whenever $(n+1)x \leq ny$ 
for $x,y \in W(A)$ and $n \in \mathbb{N}$.  Finally, let $\mathrm{sr}(A)$ denote
the stable rank of $A$.
\en

\bn
\begin{thms}\label{AHinv}
Let $A$ be a simple, unital AH algebra.  Then, $\mathrm{I}(A \otimes \mathcal{Z})$
is the augmented invariant of a simple, unital AH algebra of bounded topological dimension.
\end{thms}

\begin{nproof}
The main theorem of \cite{V3} states that every instance of the Elliott invariant
for which 
\begin{enumerate}
\item the ordered $\mathrm{K}_0$-group has the Riesz interpolation property, 
is simple and weakly unperforated, and the torsion-free part of $\mathrm{K}_{0}$ is not isomorphic to $\Z$ or the trivial group,  
\item $\mathrm{K}_1$ is a countable abelian group,
\item the tracial state space is a non-empty metrizable Choquet simplex, and
\item the continuous affine pairing between traces and states on $\mathrm{K}_0$ preserves extreme points
\end{enumerate}
occurs as $\mathrm{Ell}(B)$ for some simple, unital AH algebra $B$ of bounded topological
dimension.  Items $(2)$, $(3)$, and $(4)$ above are automatically satisfied for AH algebras (see \cite{V3}, for instance), so 
we need only establish the properties in $(1)$ for   
\[
\left((\mathrm{K}_0(A \otimes \mathcal{Z}), {\mathrm{K}_0(A \otimes \mathcal{Z})}^+, 
[\mathbf{1}_A \otimes \mathbf{1}_{\mathcal{Z}}])\right)
\]
in order to prove our theorem with $\mathrm{Ell}(A \otimes \Zh)$ in place of $\mathrm{I}(A \otimes \Zh)$. 
The full conclusion of the theorem then follows from the fact that 
the ordering on $\mathrm{K}_*(A\otimes \Zh)$ is the strict one coming from $\mathrm{K}_0(A\otimes \Zh)$ (Proposition \ref{Z-stable-invariant}).

Weak unperforation and simplicity for $\mathrm{K}_0(A \otimes \mathcal{Z})$ follow 
from Theorem 1 of \cite{GJS} and the simplicity of $A$, respectively.

The algebra $A$ is the limit of an inductive sequence $(A_j,\phi_j)$, where each
$A_j$ is a finite direct sum of $n_j$ homogeneous algebras over compact connected 
Hausdorff spaces.  For each $j \in \mathbb{N}$ and each $1 \leq i \leq n_j$, let 
\[
(G^{i,j} := G_0^{i,j} \oplus G_1^{i,j},G^{i,j+})
\]
be the graded ordered $\mathrm{K}_*$-group of the $i^{\mathrm{th}}$ building 
block of $A_j$ tensored with $\mathcal{Z}$.  By the simplicity of $A$ we may
assume, modulo compression of the inductive sequence, that the partial
morphisms
\[
\mathrm{K}_*(\phi_j)_i^k:G^{i,j} \to G^{k,j+1}
\]
have large multiplicity with respect to the co-ordinates
of $G_0^{i,j}$ and $G_0^{k,j+1}$ which correspond to the rank 
of a projection in the $i^{\mathrm{th}}$ direct summand of $A_j$
and the $k^{\mathrm{th}}$ direct summand of $A_{j+1}$.  This 
co-ordinate in each of $G_0^{i,j}$ and $G_0^{k,j+1}$ is precisely
the free co-ordinate which dominates the order in each group, as
discussed above.  Call this the first co-ordinate for convenience.
 Now let $a_1,a_2 \leq b_1,b_2$ be elements of
$G^{i,j}$, and note that this implies that either the first co-ordinate
of both $b_1$ and $b_2$ is strictly larger than either of the first co-ordinates
of $a_1$ and $a_2$, or, without loss, that $a_1 = b_1$.  In the latter case,
$a_1 = b_1$ is an interpolating element.  In the first case we push
the four elements forward via $\mathrm{K}_*(\phi_j)_i^k$, and note that
the first co-ordinate of the images of the $b_i$'s exceeds the first
co-ordinate of either of the $a_i$'s by at least two, whence, by the
discussion in \ref{torsion-free interpolation}, there is an interpolating element $g \in G^{k,j+1}$.

To see that the free part of $\mathrm{K}_0(A \otimes \mathcal{Z})$ 
cannot be cyclic, simply note that this would imply that the 
first co-ordinate multiplicities of the partial maps 
$\mathrm{K}_*(\phi_j)_i^k$ are bounded, contradicting the
simplicity of $A$.  
\end{nproof}
\en

\bn
\begin{cors}\label{weakdiv}
Let $A$ be a simple, unital, and infinite-dimensional AH algebra.  
Then, $A \otimes \Zh$ is weakly divisible in the sense of \cite{PR}.
\end{cors}

\begin{nproof}
Every simple, partially ordered abelian group $(G,G^+)$ which is weakly unperforated,
has the Riesz interpolation property, and such that $G/G_{\mathrm{tor}}$ is not cyclic must
be weakly divisible.  This statement follows from three facts:  any such group can be realised
as the ordered $\mathrm{K}_0$-group of a simple, unital AH algebra of bounded topological
dimension (\cite{V3});  any such algebra is approximately divisible (\cite{EGL2});  
approximately divisible $C^*$-algebras have weakly divisible $\mathrm{K}_0$-groups
(immediate).  It follows from Theorem \ref{AHinv} that 
\[
\left((\mathrm{K}_0(A \otimes \mathcal{Z}), {\mathrm{K}_0(A \otimes \mathcal{Z})}^+, 
[\mathbf{1}_A \otimes \mathbf{1}_{\mathcal{Z}}])\right)
\]
is weakly divisible, i.e., each $x \in {\mathrm{K}_0(A \otimes \mathcal{Z})}^+$ has
a decomposition $x=2y+3z$ for some $y,z \in {\mathrm{K}_0(A \otimes \mathcal{Z})}^+$.
The stable rank of $A \otimes \mathcal{Z}$ is one by Theorem 6.7 of \cite{R4}.  Thus,
for any projection $p \in A \otimes \mathcal{Z}$, there is a unital embedding of $\mathrm{M}_2 \oplus
\mathrm{M}_3$ into $p(A \otimes \mathcal{Z})p$.  By Lemma 5.2 of \cite{PR}, this suffices for
the weak divisibility of $A \otimes \mathcal{Z}.$ 
\end{nproof}
\en

\bn
\begin{cors}
Let $A$ be a simple, unital, and infinite-dimensional AH algebra.  If $\mathrm{K}_0(A)$
is weakly unperforated, then it is also weakly divisible and has the Riesz interpolation
property.  
\end{cors}

\begin{nproof}
As in the proof of Corollary \ref{weakdiv}, every weakly unperforated, simple, partially ordered abelian group $(G,G^+)$ with the
Riesz interpolation property is weakly divisible whenever $G/G_{\mathrm{tor}}$ is
not cyclic.  Thus, it will
suffice to prove that $(\mathrm{K}_0(A),\mathrm{K}_0(A)^+)$ has the Riesz interpolation
property, and that $\mathrm{K}_0(A)/{\mathrm{K}_{0}(A)}_{\mathrm{tor}}$ is not cyclic.

Let $A \cong \lim_{j \to \infty}(A_j, \phi_j)$ be an inductive limit decomposition of $A$,
as in the proof of Theorem \ref{AHinv}.  Notice that in this proof, the only property of
the ordered groups $(\mathrm{K}_0(A_j),\mathrm{K}_0(A_j)^+)$ required to establish the
Riesz interpolation property for $(\mathrm{K}_0(A),\mathrm{K}_0(A)^+)$ is this:   each
direct summand of $(\mathrm{K}_0(A_j),\mathrm{K}_0(A_j)^+)$ corresponding to a connected 
component of $\mathrm{Sp}(A_j)$ has the strict order coming from the cyclic subgroup generated 
by the $\mathrm{K}_0$-class of the trivial complex line bundle.  

Now suppose that $x \in (\mathrm{K}_0(A_j),\mathrm{K}_0(A_j)^+)$ has strictly positive virtual
dimension.  The stability properties of vector bundles imply that some multiple $mx$, $m \in
\mathbb{N}$, of $x$ is both positive and non-zero in $(\mathrm{K}_0(A_j),\mathrm{K}_0(A_j)^+$.
If the image of $mx$ in $(\mathrm{K}_0(A),\mathrm{K}_0(A)^+)$ is non-zero, then the image
of $x$ in the same group is positive by weak unperforation.  Thus, we may equip each direct
summand of $(\mathrm{K}_0(A_j),\mathrm{K}_0(A_j)^+)$ corresponding to a connected 
component of $\mathrm{Sp}(A_j)$ with the strict order coming from the cyclic subgroup generated 
by the $\mathrm{K}_0$-class of the trivial complex line bundle without disturbing 
$(\mathrm{K}_0(A),\mathrm{K}_0(A)^+)$.  The Riesz interpolation property now follows from
the proof of Theorem \ref{AHinv}, as does the property that 
$\mathrm{K}_0(A)/{\mathrm{K}_{0}(A)}_{\mathrm{tor}}$ is not cyclic.
\end{nproof}
\en

\bn
\label{saturation}
Given a partially ordered abelian group $(G,G^+)$, form the cone
\[
\overline{G^+} := \{x \in G \, | \, mx \in G^+ \ \mathrm{and} \ mx \neq 0 \ \mathrm{for} \ \mathrm{some} \ m \in \mathbb{N}\} \, .
\]
The order on $G$ obtained by replacing $G^+$ with $\overline{G^+}$ is called the saturated order.
It follows from Theorem 1 of \cite{GJS} that 
\begin{eqnarray*}
(\mathrm{K}_0(A \otimes \mathcal{Z}), {\mathrm{K}_0(A \otimes \mathcal{Z})}^+)  & \cong & (\mathrm{K}_0(A),\overline{\mathrm{K}_0(A)^+}) \, , 
\end{eqnarray*}
whence the saturated order on the $\mathrm{K}_0$-group
of a simple, unital, and infinite-dimensional AH algebra is weakly divisible and
has the Riesz interpolation property.
\en

\bn
\begin{cors}\label{onetracedense}
Let $A$ be a simple, unital and infinite-dimensional AH algebra with unique tracial state.  
Then, the image of $\mathrm{K}_0(A)$ in $\R = \Aff(\mathrm{T}(A))$ is dense.
\end{cors}

\begin{nproof} 
By the comments preceding the statement of the Corollary, the saturation of
$(\mathrm{K}_0(A),\mathrm{K}_0(A)^+)$ is weakly divisible, which implies that
$\mathrm{K}_0(A)$ contains elements of arbitrarily small trace.  The image
of $\mathrm{K}_0(A)$ under the unique tracial state on $A$ is therefore 
a dense subset of $\mathbb{R} \cong \Aff(\mathrm{T}(A))$, as desired.
\end{nproof}

\en

\bn
Based on Theorem 1 of \cite{Ng}, Nate Brown has shown us a proof of the following theorem, 
which will also follow from the (more general) results in \cite{Wi6}. The arguments of 
\cite{Ng} and \cite{Wi6} are quite different; the latter is built upon the methods developed in \cite{KW}, \cite{Wi4} and \cite{Wi5}. The second statement of the theorem follows from \cite{Li1}.

\begin{thms}
Let $A$ be a simple unital AH algebra. If $A\otimes \Zh$ has real rank zero (this is 
automatically the case if $A$ has real rank zero), then $A \otimes \Zh$ is tracially AF. In particular, $A$ is AH of bounded topological dimension.
\end{thms}
\en

\bn
If, in the preceding theorem, $A$ happens to have only one tracial state, then real rank zero 
follows from Corollary \ref{onetracedense} above and Corollary 7.3 of \cite{R4}:
 
\begin{thms}\label{Z-stable+uniquetrace=TAF}
Let $A$ be a simple, unital AH algebra with unique tracial state.  Then,
$A \otimes \mathcal{Z}$ is tracially AF.
\end{thms}

\begin{nproof}  Since $A \otimes \mathcal{Z}$ is $\mathcal{Z}$-stable, it has 
stable rank one by Theorem 6.7 of \cite{R4}.  As mentioned above, $A \otimes \mathcal{Z}$
has real rank zero.  Furthermore, it is locally type I with weakly unperforated
$\mathrm{K}_*$-group.  It follows from Theorem 7.1 of \cite{Br} that 
$A \otimes \mathcal{Z}$ is tracially AF.
\end{nproof}
\en

\bn The last result of this section combines results follows from \cite{Ng}, \cite{R4}, \cite{Li1}, \cite{EGL2} and our Theorem \ref{divisible-Z-stable}.  

\begin{thms}\label{classequiv}
Let $A$ be a simple, unital, and infinite-dimensional AH algebra, and suppose that
the image of $\mathrm{K}_0(A)$ in $\Aff(\mathrm{T}(A))$ is uniformly dense.  Then,
the following are equivalent:

\begin{enumerate}
\item $A$ is $\mathcal{Z}$-stable;
\item $A$ is AH of bounded topological dimension;
\item $A$ is tracially AF;
\item $W(A)$ is almost unperforated and $\mathrm{sr}(A)=1$.
\end{enumerate}
\end{thms}

\begin{nproof}
$(1) \Rightarrow (4)$.  If $A$ is $\mathcal{Z}$-stable, then $W(A)$ is almost unperforated and
$\mathrm{sr}(A)=1$ by Theorems 4.5 and 6.7 of \cite{R4}, respectively.

$(4) \Rightarrow (3)$.  $A$ is simple, unital, exact, has stable rank one, and 
$W(A)$ is almost unperforated.  Add to this the condition that
the image of $\mathrm{K}_0(A)$ in $\Aff(\mathrm{T}(A))$ is uniformly dense, and one
has the hypotheses of Proposition 7.1 of \cite{R4}.  The conclusion is that 
$A$ has real rank zero.  It then follows from Theorem 1 of \cite{Ng} that $A$
is tracially AF.

$(3) \Rightarrow (2)$.  If $A$ is tracially AF, then it is AH of bounded topological dimension by Lin's classification theorem (cf.\ \cite{Li1}, \cite{R1}).

$(2) \Rightarrow (1)$. $A$ is approximately divisible by \cite{EGL2}.
This implies $\mathcal{Z}$-stability by Theorem \ref{divisible-Z-stable}.
\end{nproof} 

Note that the class of algebras described in the hypotheses of Theorem \ref{classequiv}
need not satisfy any of the three equivalent conditions in the conclusion of the same.
Indeed, Villadsen has constructed simple, unital AH algebras of arbitrary finite
stable rank having unique trace and projections of arbitrarily small trace (\cite{V2}).

If $A$ as in the hypotheses of Theorem \ref{classequiv} has a unique tracial state, Corollary \ref{onetracedense}
implies that the density condition on the image of $\mathrm{K}_0$ can be dropped.  

There is evidence to suggest that conditions $(1)$ and $(4)$ in the conclusion of Theorem \ref{classequiv}
are equivalent for general simple, unital and infinite-dimensional AH algebras. 
\en

\section{ASH algebras}

\bn
The list of ASH classification results as mentioned in the introduction can be divided into two groups:  those results 
which cover approximately divisible $C^*$-algebras, and those which do not.  
Obviously, one may apply Theorem \ref{divisible-Z-stable} to the approximately divisible
algebras to obtain $\mathcal{Z}$-stability, so our task is twofold:  decide which
of the classification results for ASH algebras cover approximately divisible $C^*$-algebras,
and find an alternative method for proving the $\mathcal{Z}$-stability of the classified ASH
algebras which are not approximately divisible.  One can order the classification theorems covering 
ASH algebras which may fail to be approximately divisible by increasing order of 
generality:  \cite{JS1}, \cite{EJS}, \cite{T}, and \cite{M}.  We shall prove that the algebras
treated in \cite{M} are $\mathcal{Z}$-stable in Theorem \ref{AS-Z-stable} below.  As it
turns out, the remaining ASH classification results cover algebras which are approximately
divisible, though this fact is far from obvious for the algebras treated in \cite{JS2}.
\en

\bn
In \cite{R} Razak establishes the first classification result for simple, nuclear, and
stably projectionless $C^*$-algebras.  The algebras classified are simple 
inductive limits of subhomogeneous building blocks of the form
\begin{eqnarray*}
A_{n,k} & = & \mathrm{M}_n(\mathbb{C}) \otimes
\left\{ f \in \Ch([0,1],\mathrm{M}_k(\mathbb{C}))\, | \, \exists \, a \in \mathbb{C} \mbox{ such that } \right. \\
&&\left. f(0)=\mathrm{diag}(a,\ldots,a,0), \ f(1)=\mathrm{diag}(a,\ldots,a) \right\}.
\end{eqnarray*}
Let $\mathcal{R}$ denote this class of building blocks.  In \cite{Ts} Razak's 
results are generalised to cover simple inductive limits of finite direct sums 
of the building blocks above, and the range of the Elliott invariant for this 
class of algebras is computed.  
\en

\bn
The following Proposition appears to be known to a few experts, but has not
appeared in print.  

\begin{props}
Let $A = \lim_{i \to \infty}(A_i,\phi_i)$ be a simple inductive limit, where
each $A_i$ is a finite direct sum of building blocks from $\mathcal{R}$.  
Then, $A$ is approximately divisible.
\end{props}

\begin{nproof}
Let there be given a finite set $F \subseteq A$ and a tolerance $\epsilon > 0$.
We may assume without loss of generality that $F \subseteq A_1$.  We will use 
the existence and uniqueness theorems of \cite{Ts} to prove that, modulo 
compression of the inductive sequence $(A_i,\phi_i)$, there exists a map 
$\gamma:\mathrm{M}_5 \otimes A_1 \to A_2$ such that the diagram
\[
\xymatrix{
 & {\mathrm{M}_5 \otimes A_1}\ar[dr]^-{\gamma} & \\
{A_1}\ar[ur]^-{1 \otimes \mathrm{id}}\ar[rr]^-{\phi_1}&& {A_2}
}
\]
commutes up to $\epsilon$ on $F$.  The existence of $\gamma$ implies the approximate
divisibility of $A$.

The $*$-homomorphism $\mathbf{1}_{\mathrm{M}_{5}} \otimes \mathrm{id_{A_{1}}}:A_1 \to \mathrm{M}_5 \otimes A_1$
induces an isomorphism at the level of the augmented invariant $\mathrm{I}(\bullet)$.  (One
need only verify this at the level of traces, since every $A_i$ has trivial $\mathrm{K}$-groups.)
Thus, there is a map $\tilde{\gamma}:\mathrm{I}(\mathrm{M}_5 \otimes A_1) \to \mathrm{I}(A_2)$
making the diagram above commute at the level of the Elliott invariant.  The Local Existence
Theorem of \cite{Ts} implies that there exists a $*$-homomorphism
$\gamma': \mathrm{M}_5 \otimes A_1 \to A_2$ which agrees with $\tilde{\gamma}$
at the level of the Elliott invariant within a specified tolerance on a particular
finite subset of $\Aff(\mathrm{T}(A_1))$ depending only on $F$.
The Local Uniqueness Theorem of \cite{R} then implies that there is a unitary $u \in \mathcal{M}(A_2)$
such that 
\[
\gamma := u \gamma'(\,.\,) u^* : \mathrm{M}_5 \otimes A_1 \to A_2
\]
has the required property.
\end{nproof}
\en

\bn
We now recall the ASH algebras considered in \cite{M}.  Let $N, n, d_1,\ldots,d_N$ be
natural numbers such that $d_i$ divides $n$ for every $i$, and let $x_1,\ldots,x_N$ be distinct points in $\T$.  Denote by
$A(n,d_1,\ldots,d_N)$ the $C^{*}$-algebra
\[
\{f \in \Ch(\T) \otimes \mathrm{M}_n \, | \, f(x_i) \in \mathrm{M}_{d_i}, i=1,2,\ldots,N\} \, ,
\]
where $\mathrm{M}_{d_i}$ is embedded unitally in $\mathrm{M}_n$.  
Let $\mathcal{S}$ denote the collection of all such algebras.  These algebras are often 
referred to as dimension drop circles.
We will refer to the points $x_1,\ldots,x_N$ as the exceptional points of $A(n,d_1,\ldots,d_N)$.
The simple unital infinite dimensional inductive limits of finite direct sums of such
algebras are shown to be classified by the Elliott invariant in \cite{M}.
\en

\bn
\label{AS-Z-stable}
\begin{thms}
Let $A = \lim_{i \to \infty}(A_i,\gamma_i)$ be a simple unital and infinite-dimensional
inductive limit, where each $A_i$ is a finite direct sum of building blocks from $\mathcal{S}$.  
Then, $A$ is $\mathcal{Z}$-stable.
\end{thms}

\begin{nproof}
By Proposition \ref{inductive-limit-intertwining} it will suffice 
to prove the following:  given a finite set $F \subseteq A$, a dimension
drop interval $B = I[p,pq,q]$, a finite set $G \subseteq B$, and a 
tolerance $\epsilon > 0$, there is a unital embedding $\iota:B \to A$
such that
\[
||\iota(g)f-f\iota(g)|| < \epsilon, \ \ \forall f \in F, \ g \in G \, .
\]

By Lemma 9.6 of \cite{M}, we may assume that the unital $*$-homomorphism 
$\gamma_i$ is injective for every $i \in \mathbb{N}$.  We may further
assume that $F \subseteq A_1$, where
\[
A_1 = \bigoplus_{j=1}^l A_{1,j} \, ,
\]
and, for each $1 \leq j \leq l$, there are natural numbers $n_j$, $N_j$, and $d_{1,j},
\ldots,d_{1,N_j}$ such that $A_{1,j}= A(n_j,d_{1,j},\ldots,d_{N_j,j})$.

Consider the building block $A(pq,p,q)$ where the dimension drops occur at
$1,-1 \in \T$.  The fixed point algebra of $A(pq,p,q)$ under the automorphism
induced by the flip on $\T$ with fixed points $1,-1$ is isomorphic to $B$.
Assume that the exceptional points of $A_1$ are disjoint from $\{1,-1\} \subseteq \T$.
Let 
\[
\pi_B: B \longrightarrow A(pq,p,q) \otimes A_1
\]
be the $*$-monomorphism obtained by embedding $B$ into $A(pq,p,q)$ as the fixed point algebra
described above and then embedding $A(pq,p,q)$ into $A(pq,p,q) \otimes A_1$ as
$A(pq,p,q) \otimes \mathbf{1}_{A_1}$.  Let
\[
\pi_{A_1}:A_1 \longrightarrow A(pq,p,q) \otimes A_1
\]
be the embedding obtained by identifying $A_1$ with $\mathbf{1}_{A(pq,p,q)} \otimes A_1$. 
Let $\rho_j$ denote the restriction of $A(pq,p,q) \otimes A_{1,j}$ to the closed
subset 
\[
\Delta = \{(x,x) \, | \, x \in \T\} \subseteq \T^2
\]  
of its spectrum.  The image of $\rho_j$, say $D_j$, is easily seen to 
be a building block of the form 
\[
A(n_jpq,pqd_{1,j},\ldots,pqd_{N_j,j},pn_j,qn_j),
\]
which, by Corollary 3.6 and Lemma 3.9 of \cite{M}, has the same Elliott invariant 
as $A_{1,j}$.  Set
\[
D := \bigoplus_{j=1}^l D_j \, ; \ \ \rho := \bigoplus_{j=1}^l \rho_j \, .
\]

To prove the theorem it will suffice to establish the existence of a 
$*$-homomorphism $\phi:D \to A_2$ making the diagram
\[
\xymatrix{
& {D}\ar[dr]^-{\phi}& \\
{A_1}\ar[ur]^-{\rho \circ \pi_{A_1}}\ar[rr]^-{\gamma_1}&& {A_2}
}
\]
commute up to $\epsilon$ on $F \subseteq A_1$.  $B$ (and hence $G$)
can then be embedded into $D$ such that the image commutes with $\rho \circ \pi_{A_{1}}(F)$. 

The composition 
\[
\rho \circ \pi_{A_1}: A_1 \longrightarrow D
\]
induces a $\mathrm{KK}$-equivalence and an isomorphism at the level of 
the Elliott invariant.  (The details of this calculation are straightforward.
The interested reader is referred to sections 4 and 5 of \cite{M}.)  
Thus, there are a positive element $\mathbf{x}_{\phi} \in \mathrm{KK}(D,A_2)$ and an isomorphism
\[
\eta_{\phi}:\mathrm{Ell}(D) \longrightarrow \mathrm{Ell}(A_2) 
\] 
such that 
\[
[\rho \circ \pi_{A_1}] \cdot \mathbf{x}_{\phi} = [\gamma_1] \in \mathrm{KK}(A_1,A_2)
\]
and 
\[
\eta_{\phi} \circ \mathrm{Ell}(\rho \circ \pi_{A_1}) = \mathrm{Ell}(\gamma_1) \, .
\]
Let the superscript $\sharp$ denote the map induced by a $*$-homomorphism 
at the level of the Hausdorffized algebraic $\mathrm{K}_1$-group $\mathcal{U}(\bullet)/
\overline{D \mathcal{U}(\bullet)}$.  Then, one also has a morphism 
\[
\nu_{\phi}: \mathcal{U}(D)/\overline{D \mathcal{U}(D)} \to 
\mathcal{U}(A_2)/\overline{D \mathcal{U}(A_2)}
\]
such that 
\[
\nu_{\phi} \circ (\rho \circ \pi_{A_1})^{\sharp} = \gamma_1^{\sharp} \, .
\]
We have used the subscript $\phi$ above to suggest, as indeed will turn
out to be the case, that these invariant level maps above can be lifted to
a $*$-homomorphism $\tilde{\phi}:D \to A_2$.

By the simplicity of $A$, we may assume that the fibres of $A_{2}$ 
are of arbitrarily large rank at every point in the spectrum of $A_2$.  
Knowing this, we may invoke Theorem 8.1 of \cite{M}, specialising it 
to our purpose:
  
\begin{thms} (Mygind)
Let $\epsilon > 0$ and a finite set $H \subseteq \Aff(\mathrm{T}A_1)$ (to be specified) be 
given. Then, there exists a $*$-homomorphism $\tilde{\phi}:D \longrightarrow A_2$ such that
\begin{enumerate}
\item $[\tilde{\phi}] = \mathbf{x}_{\phi} \in \mathrm{KK}(D,A_2)$;
\item $\tilde{\phi}^{\sharp} = \nu_{\phi}$;
\item $||(\tilde{\phi} \circ \rho \circ \pi_{A_1}){\hat{}} \, (h)-\hat{\gamma}_1 (h)||< \epsilon||h||, \ \forall h \in H$.
\end{enumerate}
\end{thms}
\noindent

Thus, the diagram 
\[
\xymatrix{
& {D}\ar[dr]^-{\tilde{\phi}}& \\
{A_1}\ar[ur]^-{\rho \circ \pi_{A_1}}\ar[rr]^-{\gamma_1}&& {A_2}
}
\]
is at least approximately commutative (with respect to the finite set
$H \subseteq \Aff(\mathrm{T}A_1)$ of the theorem above) at the level 
of the Elliott invariant, $\mathrm{KK}$, and $\mathcal{U}(\bullet)/
\overline{D \mathcal{U}(\bullet)}$.  By enlarging $H$ if necessary,
we may apply Lemma 9.1 of \cite{M} to conclude that $\tilde{\phi}$
is injective.

To complete the proof we recall Theorem 7.7 of \cite{M}:  

\begin{thms} (Mygind)
Let $C_1$, $C_2 \in \mathcal{S}$, and let $\beta,\psi: C_1
\rightarrow C_2$ be $*$-homomorphisms inducing the same element of 
$\mathrm{KK}(C_1,C_2)$.  Let $F$ be
a finite subset of $C_1$, and let $\epsilon > 0$ be given.  There exists 
a natural number $l$ such that if $p$ and $q$ are positive integers with
$l \leq p \leq q$, if $\delta >0$, if
\begin{enumerate}
\item $\hat{\psi}(\hat{h}) > \frac{8}{p}, h \in H(C_1,l)$,
\item $\hat{\psi}(\hat{h}) > \frac{2}{q}, h \in H(C_1,p)$,
\item $||\hat{\beta}(\hat{h}) - \hat{\psi}(\hat{h})|| < \delta, h \in \tilde{H}(C_1,2q)$,
\item $\hat{\psi}(\hat{h}) > \delta, h \in H(C_1,4q)$, and
\item $D_{C_2}(\beta^{\sharp}(q^{´}_{C_1}(\nu^{C_1})),\psi^{\sharp}(q^{´}(\nu^{C_1}))) < \frac{1}{4q^2}$, 
\end{enumerate}
then there exists a unitary $u \in C_2$ such that $||\psi(f) - u^* \beta(f) 
u|| < \epsilon$ for all $f \in F$.
\end{thms}

\noindent
Many undefined objects appear in hypotheses (i)-(v).  We address this issue presently.
Replace $C_i$ with $A_i$ for $i=1,2$, $\beta$ with $\tilde{\phi} \circ \rho \circ \pi_{A_1}$, and 
$\psi$ with $\gamma_1$.  In \cite{M} it is proved that the simplicity of $A = 
\lim_{i \to \infty}(A_i,\gamma_i)$ allows one to assume that conditions (i), (ii) and (iv) are satisfied,
modulo compression of the inductive sequence for $A$.  We may assume that part (iii) of the
conclusion of Theorem 8.1 of \cite{M} holds for $H = \tilde{H}(A_1,4q)$, so that
condition (iii) is satisfied.  Condition (v) is a statement concerning the distance 
between two elements of the algebraic $\mathrm{K}_1$-group of $A_2$.  Since $D_{A_2}$ is
a metric and since 
\[
{\tilde{\phi} \circ \rho \circ \pi_{A_1}}^{\sharp}(\nu^{A_1}) = \gamma_1^{\sharp}(\nu^{A_2})
\]
 by part (ii) of the conclusion of Theorem 4.2 above,  one sees that condition (v) is satisfied.
Thus, $\gamma_1$ and $\tilde{\phi} \circ \rho \circ \pi_{A_1}$ agree up to $\epsilon$
on the finite set $F \subseteq A_1$ after conjugation by a unitary element in $A_2$.
Set $\phi := \mathrm{Ad}(u) \circ \tilde{\phi}$.  The inclusion of $B$ into $D$ composed
with $\phi$ yields a unital embedding 
\[
\iota:B \to A_2 \hookrightarrow A \, .
\]
The image $\iota(G)$ then commutes with $F$ up to $\epsilon$, as required.
\end{nproof}
\en

\bn
It remains to address the algebras classified in \cite{JS2}. 
Let $\mathcal{P}$ denote the collection of $C^*$-algebras of the form 
\begin{eqnarray*}
\lefteqn{A(a_1,\ldots,a_m;b_1,\ldots,b_l):=}\\
&&\left\{f \in \mathrm{M}_n \otimes \Ch ([0,1]) \, | \, f(0) \in \bigoplus_{i=1}^{m}
\mathrm{M}_{a_i}, f(1) \in \bigoplus_{j=1}^{l} \mathrm{M}_{b_j}\right\},
\end{eqnarray*}
where $n$, $a_1,\ldots,a_m$ and $b_1,\ldots,b_l$ are natural numbers such that
\[
\sum_{i=1}^{m} a_i = \sum_{j=1}^{l} b_j = n \, .
\]
Notice that the spectrum of such an algebra is not Hausdorff in general.  
We call the points 0 and 1 in the Hausdorffized spectrum of $A(a_1,\ldots,a_m;b_1,\ldots,b_l)$
\emph{broken endpoints}, and the spectra of the simple blocks of the fibres over 0 or 1
\emph{fractional endpoints}.  In \cite{JS2} it is shown that simple, unital, 
infinite-dimensional inductive limits of finite direct sums of members of $\mathcal{P}$ 
are classified by the Elliott invariant.
\en

\bn
\begin{thms}
Let $A = \lim_{i \to \infty}(A_i,\gamma_i)$ be simple unital and infinite-dimensional
inductive limit, where $A_i \in \mathcal{P}$, $i \in \mathbb{N}$.  Then, $A$ is 
approximately divisible.
\end{thms}

Before proceeding with the proof we recall some terminology from \cite{JS2}.  
Let $C,D$ be $C^*$-algebras.  Morphisms
\[
\kappa:\mathrm{K}_0(A) \to \mathrm{K}_0(B); \ \ \theta: \mathrm{T}(B)
\to \mathrm{T}(B)
\]
are said to be compatible if
\[
\langle e,\theta(t) \rangle = \langle \kappa(e),t \rangle
\]
for all $e \in \mathrm{K}_0(A)$ and $t \in \mathrm{T}(B)$.  A decomposition of a compatible pair
$(\kappa,\theta)$, denoted $(\kappa,\theta = \sum_j (\kappa_j,\theta_j)$, consists of
\begin{enumerate}
\item mutually orthogonal $C^*$-subalgebras $D_1,D_2,\ldots,D_n$ of $D$ such that $\mathbf{1}_D \in D_1 +
D_2 + \cdots + D_n$, and
\item a compatible pair $(\kappa_j,\theta_j)$ for $(C,D_j)$ for each $j$ satisfying
\[
\kappa = \iota_* \left(\bigoplus_{j=1}^{n} \kappa_j \right),
\]
where $\iota:D_1 + D_2 + \cdots + D_n \to D$ is the inclusion map, and 
\[
\theta(t) = \sum_{j=1}^n \theta_j(t|_{D_j}), \ \ t \in \mathrm{T}(B) \, ,
\]
where $\theta_j$ is again naturally extended.  
\end{enumerate}
Let $C = A(k_0^1,\ldots,k_0^{r_0};k_1^1,\ldots,k_1^{r_1})$.  Then, 
\[
\mathrm{K}_0(C) = \left\{(l_0^1,\ldots,l_0^{r_0},l_1^1,\ldots,l_1^{r_1}) \in \Z^s \times 
\Z^r \, | \, \sum_i l_0^i = \sum_j l_1^j \right\}.
\]

We now consider some basic compatible pairs:
\begin{enumerate}
\item Let $D = \mathrm{M}_m(\C)$, 
\[
\kappa(l_0^1,\ldots,l_0^{r_0},l_1^1,\ldots,l_1^{r_1}) = \sum_i l_0^i \, ,
\]
and $\theta(\mathrm{Tr}) = \mu$, where $\mathrm{Tr}$ is the unique normalised trace on $D$ and 
$\mu$ is any Radon probability measure on $[0,1]$.  Then, $(\kappa,\theta)$ is a compatible
pair for $(C,D)$.  Such a pair is called \emph{generic}.
\item Let $x$ be a broken endpoint of $C$.  Let $D = C(x)$,
\[
\kappa(l_0^1,\ldots,l_0^s,l_1^1,\ldots,l_1^r) = (l_0^1,\ldots,l_0^s) \, , 
\]
and $\theta(\delta_x^i)
= \delta_x^i$ for $1 \leq i \leq r_x$, where $\delta_x^i$ is the unique normalised trace
of the matrix algebra $\mathrm{M}_{k_x^i}$.  Such a pair is called \emph{broken}.
\item Let $y$ be a fractional endpoint of $C$, and let $D = C(y)$, 
\[
\kappa(l_0^1,\ldots,l_0^{r_0},l_1^1,\ldots,l_1^{r_1})= l_x^j
\]
($l_x^j$ is the co-ordinate corresponding to $y$), and $\theta(\mathrm{Tr}) = \delta_y$.  A
compatible pair equivalent to this will be called \emph{fractional} at $y$.
\end{enumerate}

\begin{nproof}
Let $A= \lim (A_i, \gamma_i)$ be a simple, unital and infinite dimensional
$C^*$-algebra, where $A_i \in \mathcal{P}$ for every $i \in \mathbb{N}$.
Let there be given a finite set $F \subseteq A$, a natural number $N$ and a tolerance $\epsilon > 0$.
We require an unital embedding of $\iota:\mathrm{M}_N \oplus \mathrm{M}_{N+1} \to A$ 
which commutes with $F$ up to $\epsilon$.  We may
assume that $F \subseteq A_1$, so that it will suffice to find an unital
embedding of $\mathrm{M}_N \oplus \mathrm{M}_{N+1}$ into  
$A_2$ which commutes with $\gamma_1(F)$ up to $\epsilon$.
In light of Lemma 2.3 of \cite{JS2}, we may assume that $A_1$ 
and $A_2$ are single building blocks.

We recall Theorem 3.7 of \cite{JS2} (local existence):

\begin{thms} (Jiang--Su)
Let $A_1$ be a splitting interval algebra.  Then, for any finite subset $F$ of 
$A_1$ and any positive tolerance $\epsilon$ there is a constant $K \in \N$
such that for any compatible pair $(\kappa,\theta)$ for $(A_1,A_2)$, where $A_2$
is a splitting interval algebra of generic fibre dimension $m$, there is a
constant $C$ depending only on $A_1$ and a $*$-homomorphism 
$\alpha:A_1 \to A_2$ of the standard form 
(cf. \cite{JS2}) which induces $\kappa$ and almost induces $\theta$ in the sense
that
\[
||\alpha^*(t)(f)-\theta(f)|| < \epsilon + \frac{CK ||f||}{m}
\]
for any $f \in F$, $t \in \mathrm{T}(A_2)$.
\end{thms}

\noindent
By Theorem 3.6 of \cite{JS2} we may assume that $K$, which depends only on
$F$ and $\epsilon$, is even.
Since we are dealing with a simple limit, we may assume that $m$ above is arbitrarily 
large.  If fact, we may assume that the dimension of each simple matrix block at the
endpoints of the spectrum of $A_2$ is arbitrarily large (cf. \cite{JS2}, Corollary 5.3).

We recall the construction of the map $\alpha$ above, and show that this map can
perturbed to a $*$-homomorphism $\phi:A_1 \to A_2$ which induces $\kappa$ and
almost induces $\theta$ as above, and whose image commutes with a unital copy of
$\mathrm{M}_N \oplus \mathrm{M}_{N+1}$ inside $A_2$.  We then apply the local uniqueness
result of \cite{JS2} to complete the proof.

For $x \in \{0,1\}$, let $Q_x:A_2 \to A_2(x)$ be the canonical evaluation map.
Applying Lemma 3.5 of \cite{JS2} to the compatible pairs
\[
(\kappa^{(x)},\theta^{(x)}) := ({(Q_x)}_* \circ \kappa, \theta \circ {(Q_x)}^*)
\]
We obtain decompositions 
\[
(\kappa^{(x)},\theta^{(x)}) = \sum_j (\kappa^{(x)}_j,\theta^{(x)}_j) \, ,
\]
where each $(\kappa^{(x)}_j,\theta^{(x)}_j)$ is either fractional, broken, or generic.  
Furthermore, these two decompositions have the same total number of summands, and the
same number of summands which are fractional pairs at the fractional endpoint $y$ of
$A_1$ for each such $y$.  Assume, therefore, that $(\kappa^{(0)}_j,\theta^{(0)}_j)$ 
is fractional at $y$ if and only if $(\kappa^{(1)}_j,\theta^{(1)}_j)$ is fractional
at $y$.  

We now group the remaining compatible pairs into ``batches'' in a manner similar to that
of the proof of Theorem 3.7 of \cite{JS2}.  In each simple block of $A_2(0)$, group the
generic pairs into batches of size $KN(2N+1)$, leaving at most $KN(2N+1)-1$ such pairs
unaligned.  Group equivalent broken pairs --- those corresponding to the same endpoint
of the spectrum of $A_1$ --- into batches of size $KN(2N+1)$, leaving at most $2KN(2N+1)-2$
unaligned pairs.  Carry out a similar grouping of the pairs in the decomposition for
$(\kappa^{(1)}_j,\theta^{(1)}_j)$.  By Lemma 5.2 of \cite{JS2} we may assume that at 
least $2/\epsilon$ batches of generic pairs occur inside each simple block of $A_2(0)$ 
and $A_2(1)$.  After batching, the total number of unaligned pairs in $A_2(0)$ does not exceed 
$(2+s_0)(KN(N+1)-1)$, where $s_0$ denotes the number of direct summands of $A_2(0)$.  
Assume that the total number of batches at 1 does not exceed that at 0, and
that if two pairs $(\kappa^{(0)}_j,\theta^{(0)}_j)$ and $(\kappa^{(0)}_l,\theta^{(0)}_l)$
are in the same batch, then so are $(\kappa^{(1)}_j,\theta^{(1)}_j)$ and 
$(\kappa^{(1)}_l,\theta^{(1)}_l)$.  

Find mutually orthogonal subalgebras $B_1, B_2, \ldots, B_k$ of $A_2$ such that:
\begin{enumerate}
\item if $(\kappa^{(0)}_j,\theta^{(0)}_j)$ is fractional at $y$, then $B_j \simeq A_1(y)$;
\item if $(\kappa^{(0)}_j,\theta^{(0)}_j)$ is not fractional, then $B_j$ is a splitting
interval algebra with $B_j(x)$ isomorphic to the fibre of $A_1$ giving rise to $\kappa^{(x)}_j$,
$x \in \{0,1\}$.
\end{enumerate}
Define $\phi$ to be equal to $\alpha$ on the batched compatible pairs for the time being.  
Note that any perturbation of $\alpha$ on the subalgebras $B_j$ corresponding to 
fractional or unaligned pairs will have a negligible effect on traces.

Consider the fibre over a  broken endpoint of the spectrum of $A_2$ having the greater number, say $\beta$, 
of simple blocks.  Choose a batch of generic pairs from each such block.  It is elementary that 
one may choose $\beta$ batches of generic pairs from among the simple blocks over the opposite
endpoint such that at least one batch is chosen from each simple block, and the proportion of 
batches chosen among all generic batches in any given simple block does not exceed $\epsilon$.  Pair these batches of
generic compatible pairs, and assume that if $(\kappa^{(0)}_j,\theta^{(0)}_j)$ is in a chosen
batch at one endpoint, then $(\kappa^{(1)}_j,\theta^{(1)}_j)$ is in the paired batch at the other
endpoint.  Let $B^{(l)}$ be the direct sum of the $B_j$ corresponding to the $l^{\mathrm{th}}$
pair of batches.  Define a $*$-homomorphism $\phi_l:A_1 \to B^{(l)}$ to be the direct sum of 
$(KN/2)(2N+1)$ copies of evaluation at each of the two broken endpoints of the spectrum of
$A_1$, $1 \leq l \leq \beta$.  Notice that this choice of morphism does not change the induced 
maps $\kappa^{(0)}_j$ and $\kappa^{(1)}_j$ corresponding to the $B_j$ under consideration.  
Let $J_{\beta}$ be the set of all indices contained in the paired batches above.   

We now address the unaligned pairs.  To each 
unaligned compatible pair in $A_2(0)$ there corresponds a compatible pair in $A_2(1)$ and a subalgebra
$B_j$ of $A_2$ whose endpoints lie entirely inside a simple block of $A_2(0)$ and $A_2(1)$, respectively.  
The map defined in Theorem 3.7 of \cite{JS2} from $A_1$ to $B_j$ can be replaced 
with evaluation at a broken endpoint, say $0$, of $A_1$ without changing the induced map on 
$\mathrm{K}$-theory.  Since the number of unaligned pairs is negligible when compared to the
total number of pairs, this modification has a negligible effect on traces.  Let $I$ denote the
set of indices of the unaligned pairs and let 
\[
\phi_{\beta+1}:A_1 \to \bigoplus_{j \in I} B_j   
\]
be the direct sum of the modified maps above.

It follows from the proof of Theorem 3.7 of \cite{JS2} that there is a morphism 
\[
\tilde{\phi}:A_1 \to \bigoplus_{j \notin J_{\beta} \cup I} B_j
\]
(which is in fact equal to $\alpha$ on these $B_j$) such that 
\[
\phi := \tilde{\phi} \oplus \left( \bigoplus_{l=1}^{\beta+1} \phi_l \right)
\]
induces $\kappa$ and almost induces $\theta$ in the manner of the conclusion of 
Theorem 3.7 of \cite{JS2}.  

Note that, by the proof of Theorem 3.7 of \cite{JS2}, the direct sum of the $B_j$ from any given 
batch sits inside a subalgebra of $A_2$ isomorphic to $\mathrm{M}_{nKN(2N+1)}(\Ch([0,1]))$, where
$n$ is the dimension of the generic fibre of $A_1$.  

The important difference between the construction of $\phi$ produced here and that of $\alpha$ in the 
proof of Theorem 3.7 of \cite{JS2} is that, using the assumption of simplicity, we are able
to guarantee that the multiplicity of the evaluation of $A_1$ at any fractional endpoint is
at least $N(N+1)$ in each simple block over the two endpoints of the Hausdorffized spectrum 
of $A_2$, and that the map $\tilde{\phi}$ can be written as the $N(2N+1)$-fold direct sum
of another morphism $\psi$.  Furthermore, the image of 
\[
 \bigoplus_{l=1}^{\beta+1} \phi_l 
\]
consists entirely of fractional evaluations.    Let $s_x$ denote the number of simple 
blocks of $A_2(x)$, $x \in \{0,1\}$.  For $1 \leq i \leq s_0$, write $a_iN+b_i(N+1)$
for the number of evaluations at $y$ in the $i^{\mathrm{th}}$ simple block of $A_2(0)$.
Since there are at least $N(N+1)$ such evaluations in each simple block, we may assume
that $a_i,b_i \geq 0$ for all $i$.  Similarly, write $c_jN+d_j(N+1)$ for the number of 
evaluations at $y$ in the $j^{\mathrm{th}}$ block of $A_2(1)$ with $c_j,d_j \geq 0$.   
Since, by the construction of the $\phi_l$, there are at least $2KN(2N+1)$ such 
evaluations in each simple block, one may modify the choices of the $a_i$, $b_i$, 
$c_j$, and $d_j$ so that
\[
\sum_i a_i = \sum_j c_j \ \ \ \mathrm{and} \ \ \ \sum_i b_i = \sum_j d_j \, .
\]
Having noticed this, one sees that $\phi$ can be factored through 
\[
A_1 \otimes (\mathrm{M}_N \oplus \mathrm{M}_{N+1}) \, .
\]  

One can easily verify the hypotheses of Theorem 4.2 of \cite{JS2} (local uniqueness)
for $\phi$, $\gamma_1$, $F$, and $\epsilon$, and conclude that $\phi$ is approximately unitarily equivalent
to $\gamma_1$ on the finite subset $F$ of $A_1$.  This establishes the approximate 
divisibility of $A$. 
\end{nproof}
\en

\end{document}